\documentclass[11pt]{article}
\usepackage{graphicx}
\usepackage{tikz}
\usepackage[pdf]{pstricks}
\usepackage{amssymb}
\usepackage{amsmath,amsthm}
\usepackage{verbatim}
\usepackage{amsfonts}
\usepackage[T1]{fontenc}
\usepackage{multicol}
\usepackage{color}
\usepackage{pdfpages}

\newtheorem{thm}{Theorem}

\newtheorem{lem}[thm]{Lemma}

\newcommand{\gs}{\mathsf{g}}

\newcommand{\Sym}{\mathcal{S}}

\newcommand{\U}{\mathcal{U}}
\newcommand{\V}{\mathcal{V}}
\newcommand{\X}{\mathcal{X}}
\newcommand{\Y}{\mathcal{Y}}

\newcommand{\Z}{\mathbb{Z}}

\title{THE HARBORTH CONSTANT OF DIHEDRAL GROUPS}
\author{NIRANJAN BALACHANDRAN\footnote{Dept. of Mathematics, Indian Institute of Technology Bombay, Mumbai, 
India. email: niranj@math.iitb.ac.in},\  ESHITA MAZUMDAR\footnote{Center for Combinatorics,
Nankai University, Tianjin 300071, P. R. China. email: eshitamazumdar@yahoo.com. Supported by NSFC, Grant No. 11681217.}, AND
KEVIN ZHAO\footnote{Center for Combinatorics,
Nankai University, Tianjin 300071, P. R. China. email: zhkw-hebei@163.com} }

\begin{document}
\maketitle
\begin{abstract} 
The Harborth constant of a finite group $G$, denoted $\gs(G)$, is  the
smallest integer $k$ such that the following holds: For  $A\subseteq G$ with $|A|=k$,  there exists  $B\subseteq A$ with $|B|=\exp(G)$ such that the elements of $B$ can be rearranged into a sequence whose product equals  $1_G$, the identity element of $G$.  The Harborth constant is a well studied  combinatorial invariant  in the case of abelian groups.  In this paper, we consider a generalization $\gs(G)$ of this combinatorial invariant for nonabelian groups and prove that if $G$ is a dihedral group of order $2n$ with $n\ge 3$, then $\gs(G) = n + 2$ if $n$ is even and $\gs(G) = 2n + 1$ otherwise.
\end{abstract}

\textbf{Keywords:}
 Zero-sum problems, Harborth constant, dihedral group. \\

2010 AMS Classification Code:  11B30, 11B75, 20D60.

\section{Introduction}For a positive integer $n$, we shall denote the set $\{1,\ldots,n\}$ by $[n]$. Let $G$ be a finite group written multiplicatively with $1=1_G$ denoting the identity element.  A sequence $\X=(x_1,\ldots,x_k)$ of elements of $G$ will be called a $G$-sequence of length $k$, or simply a sequence of length (or size) $k$, if the group in question is unambiguous. We say that a $G$-sequence $\X=(x_1,\ldots,x_k)$ is a \textit{one-product sequence} if there is a permutation $\pi$ of the indices in $[k]$ such that $\displaystyle\prod_{j=1}^{k} x_{\pi(j)} = 1$.  We say that a sequence $\X=(x_1,\ldots,x_k)$ \textit{admits a one-product subsequence} if there is a nontrivial subsequence $\Y=(x_{i_1},\ldots,x_{i_m})$ of $\X$ such that $\Y$ is a one-product sequence. If the group is abelian and is written additively, then it is customary to  refer to one-product sequences as zero-sum sequences. If $\X$ is a subset of $G$, then we may interpret $\X$ as a sequence (in no particular order), and in these instances too, we may refer to one-product subsequences of $\X$ without, hopefully, any ambiguity.

One of the most well-studied problems in combinatorial number theory is the determination of what may collectively be referred to as zero-sum constants. For instance, the Davenport constant of an abelian group $G$ is the 
minimum $k$ such that every sequence $\X$ of length $k$ admits a non-trivial zero-sum subsequence. One might also impose a size constraint on the size of the zero-subsequence that one hopes to find in a given $G$-sequence. For instance, the EGZ constant of a group $G$,  so named after a theorem of Erd\H{o}s, Ginzburg, and Ziv, is the minimum $k$ such that every $G$-sequence of length $k$ admits a zero-sum subsequence of length $\exp(G)$. The theorem of Erd\H{o}s, Ginzburg and Ziv (\cite{EGZ}) states that any $\Z_n$-sequence of length $2n-1$ admits a zero-sum subsequence of length $n$; in other words, the EGZ constant of the cyclic group $\Z_n$ is $2n-1$. 

Another set of combinatorial group invariants arises when one imposes constraints on the $G$-sequences under consideration. For instance, the Erd\H{o}s-Heilbronn problem considers determining the minimum $k$ such that every \textit{subset}  of size $k$ of an abelian group $G$  admits a non-trivial zero-sum subsequence.  Szemer\'edi (see \cite{Sze}) showed, settling a conjecture of Erd\H{o}s and Heilbronn, that there exists an absolute constant $C>0$ such that the following holds: If $G$ is an abelian group of order $n$ and $A\subseteq G$ with $|A|\ge C\sqrt{n}$ then $A$ admits a zero-sum subsequence. It is not hard to see that for the cyclic group $\Z_n$ this result is asymptotically tight: If $t(t+1)/2<n$, then the set $A=\{1,\ldots,t\}$ admits no non-trivial zero-sum subsequence.

Harborth, motivated by a problem on lattice points, (see \cite{MORS}) introduced yet another combinatorial invariant of an abelian group $G$, called its Harborth constant, which is defined as the minimum $k$ such that every subset of size $k$ of $G$ admits a zero-sum subsequence of size $\exp(G)$. However in this case, unlike the Erd\H{o}s-Heilbronn problem, it is not true that such a constant is always well defined. In fact, the group $\Z_{2n}$, admits no zero-sum subsequence of size $2n$. In such pathological cases, we adopt the local convention and define the Harborth constant to be $|G|+1$. The Harborth constant is known for a few classes of abelian groups; see \cite{MORS} for more details.
 
While all of the aforementioned invariants were defined and studied primarily for abelian groups, there has, of late been considerable interest in extending these notions over to  
non-abelian groups as well. For instance, an analogue of the Davenport constant was defined and determined for the dihedral group (see \cite{B} and
\cite{GL}). In the case of non-abelian groups one may also impose restrictions upon the structure of the permutation $\pi$ that establishes that a given sequence admits a one-product subsequence. For instance, one may require that the permutation in question is the identity;  this translates to viewing the given sequence as a word in the free monoid $G^*$ over $G$, and seeking a subword whose image under the canonical epimorphism $G^*\to G$ is the identity element of $G$ (see \cite{Loth}).  However, if we restrict our attention to subsets of $G$, then since there is no canonical sequence that represents a subset of a group, imposing a restriction on the permutation does not translate into a meaningful invariant for the group. Hence, we shall, in this paper contend with the simpler notion (as in the first paragraph of the Introduction) which we shall recall now. 

For a finite group $G$, not necessarily abelian,  the Harborth constant of $G$ is defined to be the smallest integer $k$ such that every subset  $A\subseteq G$ of cardinality $k$  admits a one-product subsequence of length $\exp(G)$, should this notion be well defined. In case it is not, then we define $\gs(G)=|G|+1$, as in the abelian case.

Our main result in this paper is the following:
\begin{thm}\label{thm1}
For any integer $n\ge 3$ and $G= D_{2n}$, 
\begin{eqnarray*} \gs(G) =\left\{    
\begin{array}{ll}
            n+2 & \textrm{if\ }2\mid n, \\
             2n+1 & \textrm{otherwise.}
         \end{array}
         \right.
\end{eqnarray*}
\end{thm}

We first prove a couple of lemmas in the next section and then prove theorem \ref{thm1} in the subsequent section. We conclude the paper with some ideas for future research.
\section{Preliminaries}
We begin by setting up some notation. We shall write the dihedral group as 
$$D_{2n}=\langle x,y \mid x^2=y^n=(xy)^2=1\rangle.$$
 For a subset $S:=\{z_1,\ldots,z_s\}\subseteq D_{2n}$, and an integer $s>0$, define 
 $$\prod\nolimits_{t}(S):=\left\{\prod_{i=1}^s z_{\pi(i)} : \pi\in\Sym(s)\right\},$$ where $\Sym(s)$ is the set of all permutations of the elements of $[s]$. For a subset $A \subseteq \Z_n$,  
  $2\cdot A$ shall denote $\{2a:a\in A\}\subset\Z_n$.

The following result is actually well known, and is quite easy. Here for a general abelian group and $A,B\subseteq G$, $A+B$ denotes the subset $\{a+b:a\in A, b\in B\}$ of $G$.
\begin{lem}\label{lm3}
Let $G$ be a finite abelian group (written additively), and let $A, B $ be
non-empty subsets of $ G$
such that $|A|+|B| \geq |G|+1.$ Then  $A+B=G.$
\end{lem}
{\bf Proof:}
For $g\in G$ and a subset $A\subseteq G$, we shall denote by $g-A$, the set $\{g-a:a\in A\}$. 

For each $g\in G$, if $|g-A|+|B| > |G|$ then
$(g-A)\cap B \neq \emptyset,$ which implies $g \in A+B$.  But since $|g-A|=|A|$, we are through.
\qed

We now turn our attention to the main results of this section - two lemmas which concern the size and structure of $\prod\nolimits_{t}(S)$ if $S$ is a subset of $D_{2n}$ of the form $\{xy^{\alpha_1},\ldots,xy^{\alpha_t}\}$. As it turns out, we need the full strength of only lemma \ref{lm2_Sodd} in the proof of theorem \ref{thm1}. But the scheme of proof of lemma \ref{lm2_Sodd} closely follows that of lemma \ref{lm2} which is easier, so it helps in understanding the arguments in lemma \ref{lm2_Sodd} better. Furthermore, as the statements of the lemmas are not too complicated, they may also be regarded as results of independent interest.

\begin{lem}\label{lm2}  Suppose $n$ is even and let $s\ge 2$.
 Let $S=\{xy^{\alpha_1}, \ldots ,xy^{\alpha_{2s}}\}$ with $\alpha_i\neq \alpha_j$ for any $1\leq i<j\leq 2s$. Then
 $$
 \left|\prod\nolimits_{2s}(S)\right|\geq s.
 $$
If equality holds, then $2s$ divides $n$ and $\{\alpha_1, \ldots ,\alpha_{2s}\}$ is a coset of the unique subgroup of $\Z_n$ of order $2s$.

\end{lem}

{\bf Proof:} We start with some terminology. Call a pair of elements $u,v\in\Z_n$ with $0\le u<v\le n-1$ a \textit{matched pair} if $2u=2v$ in $\Z_n$, or equivalently, $v=u+\frac{n}{2}$. We shall refer to $v$ as the matching mate of $u$, and vice versa. The crucial observation is this: Suppose $n$ is even, and suppose $2u=2v$ in $\Z_n$. Then either $v=u$ or $(u,v)$ form a matched pair.

 Let $S=\{xy^{\alpha_1}, \ldots ,xy^{\alpha_{2s}}\}$ with the $\alpha_i$ being pairwise distinct. Without loss of generality, we assume $0\leq \alpha_1< \alpha_2< \cdots <\alpha_{2s}\leq n-1$
Write
 \begin{equation}
 xy^{\alpha_1}xy^{\alpha_{s+1}} xy^{\alpha_2}xy^{{\alpha_{s+2}}} \cdots  xy^{\alpha_{s}}xy^{{\alpha_{2s}}}=y^{\gamma} 
 \end{equation}
so that in particular,  $\gamma=(\alpha_{s+1}+\cdots+\alpha_{2s})-(\alpha_1+\cdots+\alpha_{s})$, with this addition in $\Z_n$. By swapping $\alpha_i$ with $\alpha_{s+j}$ for $1\le i,j\le s$ in the LHS of the expression above,   it follows that 
 $$
 \{y^{\gamma}\}\cup\{y^{\gamma+2(\alpha_{i}-\alpha_{s+j})}|1\leq i,j\leq s\}\subseteq \prod\nolimits_{2s}(S).
 $$

Set \begin{eqnarray*} A&:=&\{0\}\cup\left\{\alpha_{i}-\alpha_{s+j}: 1\le i,j\le s\right\},\\
 A_0&:=&\left\{0, \alpha_{s+1}-\alpha_{s}, \alpha_{s+1}-\alpha_{s-1},\ldots, \alpha_{s+1}-\alpha_1, \alpha_{s+2}-\alpha_1,\ldots,\alpha_{2s}-\alpha_1\right\}. 
\end{eqnarray*} 
Observe that $|A_0|=2s$ and that the listing of the elements of $A_0\subseteq A$ above is in increasing order, i.e., 
$$0<\alpha_{s+1}-\alpha_{s}<\alpha_{s+1}-\alpha_{s-1}<\cdots<\alpha_{s+1}-\alpha_1<\alpha_{s+2}-\alpha_1<\cdots<\alpha_{2s}-\alpha_1. $$
The proof of the inequality of the first part of the lemma is now almost done. Indeed, by the preceding discussions,  since every $u\in \gamma+2\cdot A$ describes an element $y^{u}\in\prod_{2s}(S)$, we have   
\begin{eqnarray*}
\left|\prod\nolimits_{2s}(S)\right|&\geq& \left|\{y^{\gamma}\}\cup\{y^{\gamma+2(\alpha_{i}-\alpha_{s+j})}:1\leq i,j\leq s\}\right|\\
  &\ge&  | \gamma+2\cdot A|=|2\cdot A|\ge |2\cdot A_0|\ge s,\end{eqnarray*} where the last inequality follows from the fact that at most $2$ distinct elements of $A_0$ determine the same element of $2\cdot A_0$. This proves the first part.
  
Suppose  $\left|\prod_{2s}(S)\right|=s$.  By the preceding observations it follows that  
\begin{enumerate}
\item[(a)] $\prod\nolimits_{2s}(S) = \{y^u: u\in \gamma+2\cdot A\}$. In other words, $\prod\nolimits_{2s}(S)$ is determined by $2\cdot A$ and $\gamma$.
\item[(b)]\label{2cdotA} $|2\cdot A| = |2\cdot A_0| = s$. In particular, $2\cdot A= 2\cdot A_0$. \end{enumerate}

Consequently, for each $u\in A_0$, its matching mate is also in $A_0$. This, in particular, gives
\begin{eqnarray}\label{alpha_AP1}\alpha_{s+i} -\alpha_1=\alpha_{s+1}-\alpha_{s+2-i} +\frac{n}{2}\textrm{\ for\ all\ }1\le i\le s.\end{eqnarray}

Since $2\cdot A=2\cdot A_0$, the crucial observation made at the beginning implies that every element of $A$ either is in $A_0$ or is the matching mate of some element of $A_0$. But for each $u\in A_0$, the matching mate of $u$ is also in $A_0$, so $A=A_0$. So, for instance, since 
$$\alpha_{s+1}-\alpha_2<\alpha_{s+2}-\alpha_2<\alpha_{s+2}-\alpha_1$$  we have
\begin{eqnarray*}\alpha_{s+2}-\alpha_2=\alpha_{s+1}-\alpha_1.\end{eqnarray*} More generally, since for all $1\le i\le s-1,$
\begin{eqnarray*} \alpha_{s+i}-\alpha_2<\alpha_{s+i+1}-\alpha_2<\alpha_{s+i+1}-\alpha_1\end{eqnarray*} a simple inductive argument gives 
\begin{eqnarray}\label{alpha_AP3}\alpha_{s+i+1}-\alpha_2=\alpha_{s+i}-\alpha_1\textrm{\ for\ all\ } 1\le i\le s-1.\end{eqnarray}
  From (\ref{alpha_AP1}) and (\ref{alpha_AP3})\begin{eqnarray*}\alpha_{s+i+1} -\alpha_1&=&\alpha_{s+1}-\alpha_{s+1-i} +\frac{n}{2}\\
       &=& \alpha_{s+1}-\alpha_{s+1-i} +(\alpha_{s+1}-\alpha_1)\textrm{\ for\ }1\le i\le s-1\end{eqnarray*} so that 
 \begin{eqnarray}\label{alpha_AP4}  \alpha_{s+i+1}+\alpha_{s+1-i}=2\alpha_{s+1}\textrm{\ for\ }1\le i\le s-1.\end{eqnarray} It is now easy to see from  (\ref{alpha_AP3}), and (\ref{alpha_AP4}) that $\{\alpha_1,\ldots,\alpha_{2s}\}$ forms an arithmetic progression.  In particular, if we write $\alpha_i=a+(i-1)d$, then since these elements are pairwise distinct, we have $id\not\equiv 0\pmod n$ for any $1\leq i\leq 2s-1$. Also, since $\alpha_{s+1}-\alpha_1=\frac{n}{2}$,  it follows that $sd=\frac{n}{2}$, or equivalently, $\textrm{ord}(d)=2s$ in $\Z_n$. This completes the proof.\qed

The next lemma considers the case  where $n$ is even, and $|S|$ is odd.  

\begin{lem}\label{lm2_Sodd}
Suppose $n$ is even and let $S=\{xy^{\alpha_1}, \ldots ,xy^{\alpha_{2s+1}}\}$ with $\alpha_i\neq \alpha_j$ for any $1\leq i<j\leq 2s+1$. Then
 $$
 \left|\prod\nolimits_{2s+1}(S)\right|\geq s+1.
 $$
If equality holds then $2s+2$ divides $n$ and there is a coset $K$ of the subgroup $H$ of $\Z_n$ of order $2s+2$ such that $\{\alpha_1, \ldots ,\alpha_{2s+1}\}\subset K.$ 
\end{lem}
\textbf{Remark:} The statement of lemma \ref{lm2_Sodd}, basically states that $S$ is `almost' a coset, i.e., it misses exactly one element of a coset $K$ of the subgroup $H$ of order $2s+2$. In fact, the missing element could be any one of the $2s+2$ elements of $K$, as the proof will show. 

{\bf Proof:}  As in the proof of lemma \ref{lm2}, let us assume without loss of generality that $0\le\alpha_1<\alpha_2<\cdots<\alpha_{2s+1}\le n-1$. Write
 \begin{equation}
 xy^{\alpha_1}xy^{\alpha_{s+2}} xy^{\alpha_2}xy^{{\alpha_{s+3}}} \cdots  xy^{\alpha_{s}}xy^{\alpha_{2s+1}}xy^{\alpha_{s+1}}=xy^{\alpha}\in\prod\nolimits_{2s+1}(S)
 \end{equation}
so that $\alpha=(\alpha_1+\cdots+\alpha_{s+1})-(\alpha_{s+2}+\cdots+\alpha_{2s+1})$. Again, as in the proof of lemma \ref{lm2}, we may swap $\alpha_j$ and $\alpha_{s+1+i}$, to get 
$$xy^{\alpha+2(\alpha_{s+1+i}-\alpha_j)}\in\prod\nolimits_{2s+1}(S)\textrm{\ for\ any\  } 1\le i\le s, 1\le j\le s+1.$$  
We again set
\begin{eqnarray*} A&=&\{0\}\cup\{(\alpha_{s+1+i}-\alpha_j): 1\le i\le s, 1\le j\le s+1\},\\
 A_0&:=&\{0, \alpha_{s+2}-\alpha_{s+1}, \alpha_{s+2}-\alpha_{s}, \ldots, \alpha_{s+2}-\alpha_1, \alpha_{s+3}-\alpha_1,\ldots,\alpha_{2s+1}-\alpha_1\}.\end{eqnarray*}
 and note that 
\begin{eqnarray*} \left|\prod\nolimits_{2s+1}(S)\right| &\ge & |2\cdot A|,\\
            |A_0| &=& 2s+1.\end{eqnarray*}
 As in the proof of lemma \ref{lm2}, since
 \begin{eqnarray}\label{orderS} 0<\alpha_{s+2}-\alpha_{s+1}<\alpha_{s+2}-\alpha_{s}<\cdots<\alpha_{s+2}-\alpha_1<\alpha_{s+3}-\alpha_1<\cdots<\alpha_{2s+1}-\alpha_1<n
\end{eqnarray} 
 it follows that \begin{eqnarray}\label{sigmaSodd}\left|\prod\nolimits_{2s+1}(S)\right|\ge |2\cdot A_0|\ge s+1\end{eqnarray} where the last inequality follows since $|A_0|=2s+1$. This yields the stated inequality. 
 
  Suppose equality holds in  (\ref{sigmaSodd}).  Then, exactly as in the proof of lemma \ref{lm2}, $\prod\nolimits_{2s+1}(S)$ is determined by $\alpha$ and $2\cdot A$, and $2\cdot A=2\cdot A_0$. Since  $|A_0|=2s+1$, it follows that there are exactly $s$ matched pairs among the elements of $A_0$ which leaves exactly one element of $A_0$ whose matching mate is not in $A_0$; we shall call this, a \textit{distinguished} element. As before, since $2\cdot A=2\cdot A_0$, we have 
  $$A_0\subseteq A\subseteq B_0:=A_0\cup\{\textrm{matching\ mate\ of\ the\ distinguished\ element\ of\ }A_0\}.$$
Since $A_0$ has a unique distinguished element, we proceed by fixing a possibility for the distinguishing element and exploring the consequences. As it turns out, every element of $A_0$ is a plausible choice for being the distinguished element; in fact, knowing the distinguished element determines $A$ in a unique sense, as we shall see. 

To keep the description convenient, we introduce some further terminology.  We shall refer to elements of the form $\alpha_{s+2}-\alpha_{s+2-i}$ (with $0\le i\le s$) as \textit{former} elements, and elements of the form $\alpha_{s+2+i}-\alpha_1$ (with $0\le i\le s-1$) as the \textit{latter} elements of $A_0$.  We shall call $\alpha_{s+2}-\alpha_{s+2-i}$ for $2\le i\le s-2$ as \textit{generic former elements}, and  the remaining former elements will be called \textit{special former elements}.  We shall consider here in detail, the case when $\alpha_{s+2}-\alpha_{s+2-i}$ is the distinguished element for a generic former element of $A_0$. The proofs of the other cases are very similar, and we shall relegate those details to the Appendix.

We shall again attempt to sandwich elements of $A$ between consecutive elements of $A_0$, and then reach an identity as we did in the proof of lemma \ref{lm2}. But what gets a little complicated here is that unlike in the proof of lemma \ref{lm2}, the identities extend only till the `barrier', viz., till we encounter the predecessor of the matching mate of the distinguished element (See figure 1 for an illustration). If the elements of $A$ being sandwiched move clockwise towards zero among the sequence of identities we establish, we call such a sequence of identities, a `backward propagation', and if the elements of $A$ sandwiched move closer to $\alpha_{2s+1}-\alpha_1$, we shall refer to those identities as a `forward propagation'. This terminology will help us describe our proof better. 

Let us get into the details now. First, suppose that for some $1\le i\le s$, $\alpha_{s+2}-\alpha_{s+2-i}$ is the distinguished element.  Then by the observations made earlier, and following the same kind of argument as in the proof of lemma \ref{lm2} we conclude that  
 \begin{eqnarray}\label{alpha_1}
 \alpha_{s+2}-\alpha_{s+2-j}+\frac{n}{2} = \left\{    
\begin{array}{ll}
            \alpha_{s+1+j}-\alpha_1 \hspace{0.2cm}\textrm{\ if\ }& j>i,\\ 
             \alpha_{s+2+j}-\alpha_1  \hspace{0.2cm}\textrm{\ if\ }& 0\le j<i.
         \end{array}
         \right.
\end{eqnarray}
Using $\alpha_{s+2}-\alpha_1=\frac{n}{2}$ in conjunction with  (\ref{alpha_1}) gives
 \begin{eqnarray}\label{alpha_3}
 \alpha_{s+2}-\alpha_{s+2-j} = \left\{    
\begin{array}{ll}
            \alpha_{s+1+j}-\alpha_{s+2} \hspace{0.2cm}\textrm{\ if\ }& j>i,\\ 
             \alpha_{s+2+j}-\alpha_{s+2}  \hspace{0.2cm}\textrm{\ if\ }& 0\le j<i.
         \end{array}
         \right.
\end{eqnarray}

Now, suppose specifically that $2\le i\le s-2$.  Again,  as in the proof of lemma \ref{lm2}, since $\alpha_{s+2}-\alpha_2<\alpha_{s+3}-\alpha_2<\alpha_{s+3}-\alpha_1$, the element $\alpha_{s+3}-\alpha_2\in B_0$, so $\alpha_{s+3}-\alpha_2=\alpha_{s+2}-\alpha_1$.  More generally, a simple inductive argument gives
\begin{eqnarray}\label{alpha_2a}
\alpha_{s+2+j}-\alpha_2&=&\alpha_{s+1+j}-\alpha_1\hspace{0.5cm}\textrm{\ for\ }1\le j\le i-1, \hspace{0.3cm}(\textrm{through\ forward\ propagation})\\
\label{alpha_2b} \alpha_{s+3}-\alpha_{s+3-\ell}&=&\alpha_{s+2}-\alpha_{s+2-\ell}\hspace{0.1cm}\textrm{\ for\ } 1\le\ell\le s+1\hspace{0.4cm}(\textrm{through\ backward\ propagation}).
\end{eqnarray}

 Hence (\ref{alpha_2a}), (\ref{alpha_2b}), give 
$$\alpha_{s+i+1}-\alpha_{s+i}=\cdots=\alpha_{s+3}-\alpha_{s+2}=\alpha_{s+2}-\alpha_{s+1}=\cdots=\alpha_{2}-\alpha_{1}=d\textrm{\ (say)}.$$ Also, since $\alpha_{s+2}-\alpha_1=\frac{n}{2}$, we have $d=\frac{n}{2(s+1)}$.

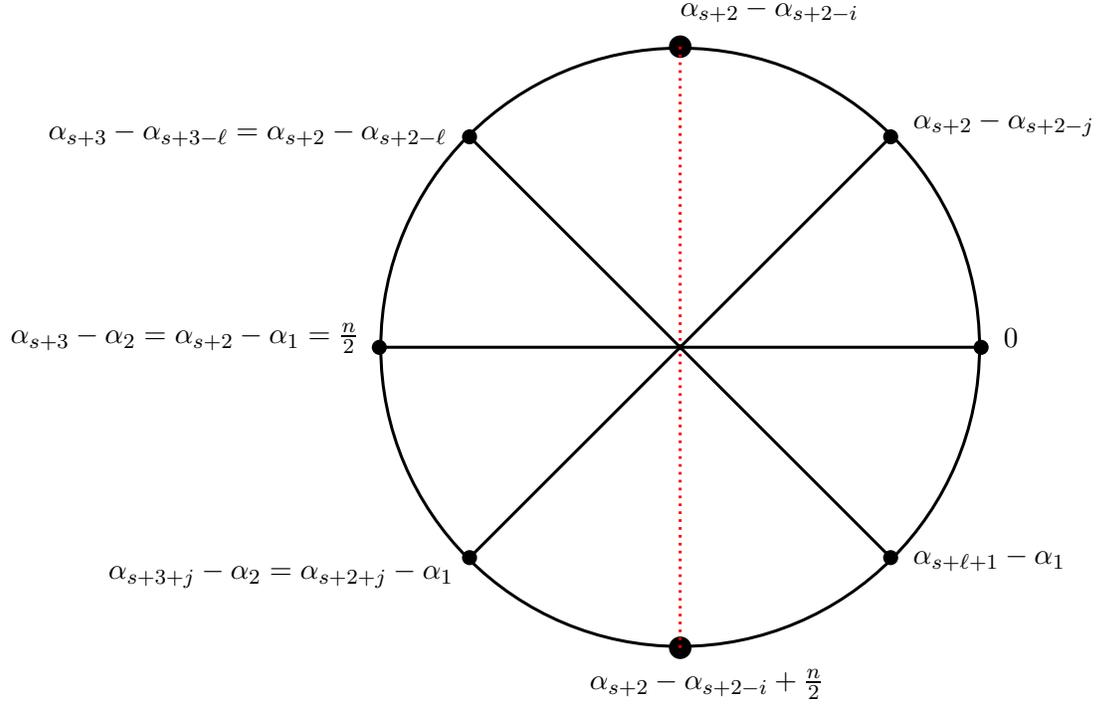
\begin{figure}
\begin{pspicture}(0,-4.72)(14.99,4.72)
\pscircle[linecolor=black, linewidth=0.04, dimen=outer](9.6,0.08){4.0}
\pscircle[linecolor=black, linewidth=0.04, dimen=outer](9.6,0.08){0.04}
\psdots[linecolor=black, dotsize=0.2](13.6,0.08)
\psdots[linecolor=black, dotsize=0.2](5.6,0.08)
\psdots[linecolor=black, dotsize=0.2](12.4,2.88)
\psdots[linecolor=black, dotsize=0.2](6.8,-2.72)
\psdots[linecolor=black, dotsize=0.2](6.8,2.88)
\psdots[linecolor=black, dotsize=0.2](12.4,-2.72)
\psdots[linecolor=black, dotsize=0.3](9.6,4.08)
\psdots[linecolor=black, dotsize=0.3](9.6,-3.92)
\psline[linecolor=black, linewidth=0.04](13.6,0.08)(5.6,0.08)
\psline[linecolor=black, linewidth=0.04, linestyle=dashed, dash=1pt 2pt, linecolor=red](9.6,4.08)(9.6,-3.92)
\psline[linecolor=black, linewidth=0.04](12.4,2.88)(6.8,-2.72)
\psline[linecolor=black, linewidth=0.04](6.8,2.88)(12.4,-2.72)
\rput[bl](13.9,0.08){$0$}
\rput[bl](12.7,2.88){$\alpha_{s+2}-\alpha_{s+2-j}$}
\rput[bl](9.6,4.4){$\alpha_{s+2}-\alpha_{s+2-i}$}
\rput[bl](1.2,2.75){$\alpha_{s+3}-\alpha_{s+3-\ell}=\alpha_{s+2}-\alpha_{s+2-\ell}$}
\rput[bl](0.7,0.0){$\alpha_{s+3}-\alpha_2=\alpha_{s+2}-\alpha_1=\frac{n}{2}$}
\rput[bl](2.0,-3.12){$\alpha_{s+3+j}-\alpha_2=\alpha_{s+2+j}-\alpha_1$}
\rput[bl](8.4,-4.6){$\alpha_{s+2}-\alpha_{s+2-i}+\frac{n}{2}$}
\rput[bl](12.7,-2.9){$\alpha_{s+\ell+1}-\alpha_1$}
\end{pspicture}
\caption{An illustration for lemma \ref{lm2_Sodd} representing the elements of $\Z_n$ cyclically. The distinguished element is $\alpha_{s+2}-\alpha_{s+2-i}$. The matched pairs among the elements of $A_0$ are joined by lines, while the distinguished element and its matching pair are joined by a dashed line. The matching mate of the distinguished element is not in $A_0$.}
\end{figure}

Since $$\alpha_{s+i}-\alpha_1=\alpha_{s+i+1}-\alpha_2<\alpha_{s+i+2}-\alpha_2<\alpha_{s+i+2}-\alpha_1,$$ it follows that 
$$\alpha_{s+i+2}-\alpha_2\in\left\{\alpha_{s+i+1}-\alpha_1, \alpha_{s+2}-\alpha_{s+2-i}+\frac{n}{2}\right\}.$$ We are precisely at the barrier that we mentioned earlier.

We now claim that $\alpha_{s+i+2}-\alpha_2=\alpha_{s+2}-\alpha_{s+2-i}+\frac{n}{2}$. Indeed, suppose if possible, that
$\alpha_{s+i+1}-\alpha_2=\alpha_{s+i+1}-\alpha_1$. Then $\alpha_{s+i+2}-\alpha_{s+i+1}=\alpha_{2}-\alpha_1=d$, so that 
$$(s+i+1)d = \alpha_{s+i+2}-\alpha_1 > \alpha_{s+2}-\alpha_{s+2-i}+\frac{n}{2}=(i+s+1)d$$ and that is a contradiction. The strict inequality above follows from the facts that $\alpha_{s+2}-\alpha_{s+2-i}+\frac{n}{2}$ is the matching mate of the distinguished term,  $ \alpha_{s+i+2}-\alpha_1$ is the matching mate of $\alpha_{s+2}-\alpha_{s+2-(i+1)}$, and  $d=\frac{n}{2(s+1)}$.

\begin{figure}
\begin{pspicture}(0,-0.32)(17.417116,0.32)
\psline[linecolor=black, linewidth=0.04](0.19711533,0.08)(16.597115,0.08)(16.597115,0.08)
\psdots[linecolor=black, dotsize=0.2](0.19711533,0.08)
\psdots[linecolor=black, dotsize=0.2](0.9971153,0.08)
\psdots[linecolor=black, dotsize=0.2](6.1971154,0.08)
\psdots[linecolor=black, dotsize=0.2](6.997115,0.08)
\psdots[linecolor=black, dotsize=0.2](7.7971153,0.08)
\psdots[linecolor=black, dotsize=0.2](10.5971155,0.08)
\psdots[linecolor=black, dotsize=0.2](11.397116,0.08)
\psdots[linecolor=black, dotsize=0.2](12.997115,0.08)
\psdots[linecolor=black, dotsize=0.2](15.797115,0.08)
\psdots[linecolor=black, dotsize=0.2](16.597115,0.08)
\rput[bl](0.19711533,-0.4){$\alpha_1$}
\rput[bl](0.9971153,-0.4){$\alpha_2$}
\rput[bl](6.1971154,-0.4){$\alpha_{s+1}$}
\rput[bl](6.997115,-0.4){$\alpha_{s+2}$}
\rput[bl](7.7971153,-0.4){$\alpha_{s+3}$}
\rput[bl](10.5971155,-0.4){$\alpha_{s+i}$}
\rput[bl](11.397116,-0.4){$\alpha_{s+i+1}$}
\rput[bl](12.997115,-0.4){$\alpha_{s+i+2}$}
\rput[bl](15.797115,-0.4){$\alpha_{2s}$}
\rput[bl](16.597115,-0.4){$\alpha_{2s+1}$}
\rput[bl](0.59711534,0.2){$d$}
\rput[bl](6.5971155,0.2){$d$}
\rput[bl](7.397115,0.2){$d$}
\rput[bl](10.997115,0.2){$d$}
\rput[bl](12.197115,0.2){$2d$}
\rput[bl](16.197115,0.2){$d$}
\end{pspicture}
\caption{The $\alpha_i$ in increasing order. If $H$ denotes the subgroup of $\Z_n$ of order $2(s+1)$, then the set $\{\alpha_1,\ldots,\alpha_{2s+1}\}$ is contained in some coset $K$ of $H$. The only missing element from $K$ is the one corresponding to the element $(s+i+1)d$ in $H$.}
\end{figure}
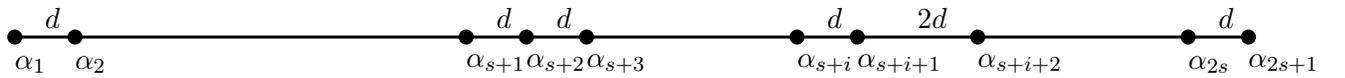

Now that we have established $\alpha_{s+i+2}-\alpha_2=\alpha_{s+2}-\alpha_{s+2-i}+\frac{n}{2}$, the remainder of the proof again follows the same line of argument that was outlined at the beginning. Since
$$ \alpha_{s+i+2}-\alpha_2<\alpha_{s+i+3}-\alpha_2<\alpha_{s+i+3}-\alpha_1$$we have (by forward propagation) $\alpha_{s+i+3}-\alpha_2=\alpha_{s+i+2}-\alpha_1$, and more generally, 
$$ \alpha_{s+\ell+1}-\alpha_2=\alpha_{s+\ell}-\alpha_1 \textrm{\ for\  all\ } i+2\le\ell\le s.$$ 
Consequently, we have the situation as in figure 2. Indeed, it is now easy to see that $\alpha_{j+1}-\alpha_j=d$ for all $j\ne s+i+1$, and since $\alpha_{s+i+2}-\alpha_2=\alpha_{s+2}-\alpha_{s+2-i}+\frac{n}{2}$, we have
\begin{eqnarray*} (s+i+1)d=\alpha_{s+2}-\alpha_{s+2-i}+\frac{n}{2}&=&\alpha_{s+i+2}-\alpha_2\\
&=&\alpha_{s+i+2}-\alpha_{s+i+1} +\alpha_{s+i+1}-\alpha_2\\ 
&=&\alpha_{s+i+2}-\alpha_{s+i+1} + (s-1)d\end{eqnarray*} which gives $\alpha_{s+i+2}-\alpha_{s+i+1}=2d$, and proves the lemma in this case.

As mentioned before, we relegate the discussion of the remaining details  of the proof to the Appendix.
\qed
\section{Proof of theorem \ref{thm1}}
Recall that $D_{2n}=\langle x,y\mid x^2=y^n=(xy)^2=1\rangle$ as mentioned in the preceding section. Let $H$ denote the cyclic subgroup of $D_{2n}$ generated by $y$ and set $N= D_{2n}\setminus H.$\\
We start with the easier case, viz., the case where $n$ is odd. Note that in this case $\exp(D_{2n})= 2n.$ From the definition of the Harborth constant, it follows that $\gs(D_{2n})\geq 2n.$ But note that for any sequence $(xy^{\alpha_1},\ldots,xy^{\alpha_t})$ the product $$(xy^{\alpha_1})\cdots(xy^{\alpha_t})=xy^{\alpha}$$ for some $\alpha\in\Z_n$ if $t$ is odd. In particular, it follows that  for every ordering $(g_1,\ldots,g_{2n})$ of the elements of $D_{2n}$, the product $\displaystyle\prod_{i=1}^{2n} g_i\ne 1$, since there are an odd number of elements of the form $xy^{\alpha}$.  Hence $\gs(D_{2n})= 2n+1$ when $n$ is odd.

We now turn to the non-trivial part of the theorem. So, in the rest of the proof, we assume that $n$ is even. 

 Consider the set $A=H\cup \{x\}$ with $H\subseteq D_{2n}$ denoting the cyclic subgroup of order $n$. Since $\displaystyle\prod_{i=0}^{n-1} y^i\ne 1$ it easily follows that for any sequence $(g_1,\ldots,g_n)$ of distinct elements from $A$, $\prod _i g_i\ne 1$, so this proves that $\gs(D_{2n})\geq n+2,\text{ if } 2\mid n.$ So in order to complete the proof, it suffices to show that 
$\gs(D_{2n})\leq n+2$, if  $2\mid n.$ 

Towards that end, suppose $S=\{xy^{u_1},\ldots,xy^{u_t},y^{v_1},\ldots,y^{v_s}\}\subseteq D_{2n}$ with $|S|=t+s=n+2$. Here, the $u_i$ (resp $v_j$) are pairwise distinct and are elements of $\Z_n$, so we may write $0\le u_i,v_j\le n-1$. Set $\U=\{u_1,\ldots,u_t\}$ and $\V=\{v_1,\ldots, v_s\}$. Note that these are viewed as subsets of $\Z_n$. Let us also for brevity's sake, write  $\sigma(S):=\prod\nolimits_{n+2}(S)$. 
We have the following subcases:
\begin{itemize}
\item {\bf $t$ is odd:} For starters, note that $t\ge 3$. Let $0\leq \gamma \leq n-1$ be such that $$\displaystyle\prod_{i=1}^{t}xy^{u_i}\prod_{j=1}^{s}y^{v_j}= xy^{\gamma}\in\sigma(S).$$ Now reorder the $u_i$ and $v_j$ as follows: $u_1'=u_2, u_2'= u_1$, and $u_i'=u_i$ for $i\ge 3$, and make a similar definition of $v_j'$ from the $v_j$, and set $\displaystyle\prod_{i=1}^{t}xy^{u_i'}\prod_{j=1}^{s}y^{v_j'}=xy^{\delta}\in\sigma(S)$. These two distinguished elements will play a pivotal role in our proof.

The starting observation is that every element in $\sigma(S)$ is of the form $xy^{\beta}$ for some $\beta\in\Z_n$, and that all these $\beta$ (i.e., when  $xy^{\beta}\in\sigma(S)$) have the same parity. 

\vspace{0.2cm}
 
 Suppose  $U, V\subseteq \Z_n$ satisfy $|U|+|V|=n+1$ (in particular, 
both $U,V\neq\emptyset$). Then by lemma \ref{lm3}, for any $\beta\in \Z_n$, there exist $u\in U,v\in V$ such that $u+v=\beta$.  As a consequence we see: For any $xy^{\beta}\in \sigma(S)$ with $1\leq \beta\leq n-1$  there exist $(u_1(\beta), v_1(\beta))$ and $(u_2(\beta),v_2(\beta))$ with the $u_i$'s and the $v_i$'s distinct, such that $u_i(\beta)+v_i(\beta)=\beta$ for $i=1,2$. We shall simply write $u_1, v_1$ (resp. $u_2,v_2$) instead of $u_i(\beta),v_i(\beta)$ for simplicity, when the $\beta$ in question is clear from the context. 

Fix $\beta$ with $xy^{\beta}\in \sigma(S)$, and suppose $(u_1,v_1), (u_2,v_2)$ are the pairs in $\Z_n$ such that $u_i+v_i=\beta$, for $i=1,2$.  Let $a_1,\ldots,a_{\ell}\in\Z_n$ be such that 
$$\{y^{a_1},\ldots,y^{a_{\ell}}\}=\prod\nolimits_{t-1}(xy^{u_2},\ldots,xy^{u_t})$$  
$\textrm{\ and\ let\ } y^{a_1}=\displaystyle\prod_{i=2}^{t}xy^{u_i}$. Again, note that all the $a_i$ have the same parity. 

Since $2\mid (t-1)$,  lemma \ref{lm2} implies that $\ell\geq \frac{t-1}{2}.$ Also, by the definition of $\gamma$, it follows that 
\begin{eqnarray*}
xy^{\gamma} &=& xy^{u_1}\prod_{i=2}^{t}xy^{u_i}y^{v_1}\prod_{j=2}^{s}y^{v_j}=xy^{u_1}\left(\prod_{i=2}^{t-1}xy^{u_i}\right)y^{-v_1}xy^{v_t}\left(\prod_{j=2}^{s}y^{v_j}\right)
  =\cdots \\ &=&  xy^{u_1}y^{v_1} \prod_{i=2}^{t}xy^{u_i}\prod_{j=2}^{s}y^{v_j}=xy^{u_1}y^{v_1}y^{a_1}\prod_{j=2}^s y^{v_j}. \end{eqnarray*} Again, this follows because  $2 \mid (t-1).$ Therefore,  
$$xy^{u_1}y^{v_1}y^{a_i}\prod_{j=2}^{s}y^{v_j}=xy^{\gamma-a_1+a_i} \in \sigma (S)
\text{ for } 1\leq i\leq \ell.$$
Similarly, 
\begin{eqnarray*}xy^{u_1}y^{v_1}\left(\prod_{i=2}^{t-1}xy^{u_i}\right)y^{v_k}xy^{u_t}\left(\prod_{\substack{2\le i\le s\\ i\neq k}}y^{v_i}\right) &=& xy^{u_1}y^{v_1}\left(\prod_{i=2}^{t-1}xy^{u_i}\right)xy^{u_t}y^{-v_k}\left(\prod_{\substack{2\le i\le s\\ i\neq k}}y^{v_i}\right)\\ 
&=&  xy^{u_1}y^{v_1} \left(\prod_{i=2}^{t}xy^{u_i}\right)\left(\prod_{j=2}^{s}y^{v_j}\right) y^{-2v_k} \\
&=& xy^{\gamma-a_1+a_i-2v_k} \in \sigma (S)\end{eqnarray*}
for all $2\leq k\leq s.$\\

Consider the set $\{-2v_2,\ldots, -2v_s\}\subseteq \Z_n.$ Since $n$ is even,  by arguments that have appeared before, it follows that  $$\Big|\{-2v_2, \ldots, -2v_s\}\Big| \geq \frac{s-1}{2}.$$

Let $A=A_{\beta,\gamma}:=\{\gamma-a_1+a_i-\beta : 1\leq i\leq \ell\}$ and $B=\{-2v_2,\ldots, -2v_s\}$. Note that both $A,B\subseteq G:=2\cdot \Z_{n}\subseteq\Z_n$ and  $$|A|+|B| \geq \ell+ \frac{s-1}{2}\geq \frac{t-1}{2}+\frac{s-1}{2}= \frac{n}{2}.$$

If $|A|+|B| > \frac{n}{2}$ then by lemma \ref{lm3} it follows that there exist $i, k$ such that $\gamma-a_1+a_i -2v_k - \beta=0$ in $\Z_n$, so 
\begin{eqnarray} \label{1prod} xy^{u_1}y^{v_1}\left(\prod_{i=2}^{t-1}xy^{u_i}\right)y^{v_k}xy^{u_t} \left(\prod_{\substack{2\le i\le s\\ i\neq k}}y^{v_i}\right) =xy^{\gamma-a_1+a_i-2v_k} = xy^{\beta}.\end{eqnarray}
But since $u_1+v_1 = \beta $, the product of the first two elements in the LHS of (\ref{1prod}) is precisely $xy^{\beta}$, and that gives us the one-product of size $n$. We are done in this scenario.


In a similar vein, let
$$\{y^{b_1},\ldots,y^{b_m}\}:=\displaystyle\prod\nolimits_{t-1}(xy^{u_2'},xy^{u_3'},\ldots,xy^{u_t'}).$$ Then proceeding exactly as before, setting $A'=A_{\beta,\delta}:=\{\delta-a_1+a_i-\beta : 1\leq i\leq \ell\}$ and $B'=\{-2v_1,-2v_3,\ldots, -2v_s\}$ we again have $|A'|+|B'|\ge\frac{n}{2}$. Thus, the only case when the proof is not complete is when 
\begin{eqnarray*}
|A|&=&\frac{t-1}{2}\hspace{0.2cm}=\hspace{0.2cm}\left|\prod\nolimits_{t-1}(xy^{u_2},\ldots,xy^{u_t})\right|,\\ |B|&=&\frac{s-1}{2}\hspace{0.2cm}=\hspace{0.2cm}|\{-2v_2,\ldots, -2v_s\}|,\\
|A'|&=&\frac{t-1}{2}\hspace{0.2cm}=\hspace{0.2cm}\left|\prod\nolimits_{t-1}(xy^{u_2'},\ldots,xy^{u_t'})\right|,\\ |B'|&=&\frac{s-1}{2}\hspace{0.2cm}=\hspace{0.2cm}|\{-2v_1,-2v_3,\ldots, -2v_s\}|.\end{eqnarray*}  

Under this scenario, it follows that both $\V\setminus\{v_1\}$ and $\V\setminus\{v_2\}$ consist of $\frac{s-1}{2}$ matched pairs from $\Z_n$. But then it is easy to see that this implies that $v_1=v_2$, and that is a contradiction.

\item {\bf $t$ is even:}  
 For starters, observe that for any permutation $\pi$ of the elements of $[t]$, we have 
$$\displaystyle\prod_{i=1}^{t}xy^{u_{\pi(i)}}\prod_{j=1}^{s}y^{v_j}=y^{\gamma(\pi)}$$ for some $\gamma(\pi)\in\Z_n$.  Moreover, all the $\gamma(\pi)$ have the same parity. Fix $0\le \gamma\le n-1$ so that \begin{eqnarray}\label{defn_gamma}\displaystyle\prod_{i=1}^{t}xy^{u_i}\prod_{j=1}^{s}y^{v_j}=y^{\gamma}.\end{eqnarray}

Since $t+s=n+2$, $\max\{t,s\} \geq \frac{n}{2}+1.$

If $s \geq \frac{n}{2}+1$ set
\begin{eqnarray*}A=\left\{    
\begin{array}{ll}
            \V\setminus\{\gamma/2\},\textrm{\ \ if\   } \gamma\in 2\cdot \Z_n,\\ \V,\hspace{1.5cm}\textrm{\ otherwise,}
         \end{array}\right. 
         \end{eqnarray*}  
         and 
$B=\{v_1,\ldots,v_s\}$.  Again, by lemma \ref{lm3}, it is easy to see that $A+B = \Z_n$, so in particular, there exist $i_0\neq j_0$ (since $\gamma/2\not\in A$ by the definition of $A$) such that $v_{i_0} + v_{j_0} = \gamma$.
Cancelling $y^{v_{i_0}}, y^{v_{j_0}}$ from among the terms in $\displaystyle\prod_{j=1}^{s}y^{v_j}$ in the LHS gives us a one-product of size $n$, and we are done in this case.

\vspace{0.2cm}

Suppose then that $s \le\frac{n}{2}$, so in particular $t \ge \frac{n}{2}+2$. In particular, for both $i=0,1$, there is a pair $(u_i,u_i')$ among $\{u_1,\ldots,u_t\}$ such that $u_i-u_i'\equiv i\pmod 2$. Hence we may without loss of generality assume that $u_2-u_1\equiv\gamma\pmod 2$ where $\gamma$ is as in (\ref{defn_gamma}). Let $$\{xy^{a_1},\ldots,xy^{a_{\ell}}\}:=\displaystyle\prod\nolimits_{t-3}(xy^{u_3},\ldots, xy^{u_{t-1}}).$$ Again, it is easy to check that all the $a_i$ have the same parity.  Set $$xy^{a_1}=\prod_{i=3}^{t-1} xy^{u_i}.$$

Now, for any $1\le i\le\ell$,
\begin{eqnarray*} xy^{u_1}xy^{u_2}xy^{a_i}y^{v_k}(xy^{u_t})\prod_{\substack {1\le j\le s\\j\ne k}}y^{v_j}&=& y^{u_2-u_1}xy^{a_1}y^{a_i-a_1+v_k}(xy^{u_t})\prod_{\substack {1\le j\le s\\j\ne k}}y^{v_j}\\ &=& y^{u_2-u_1}xy^{a_1}(xy^{u_t})\bigg(\prod_{1\le j\le s}y^{v_j}\bigg)y^{-a_i+a_1-2v_k}\\&=& y^{\gamma+a_1-a_i-2v_k}\in \sigma(S),
\end{eqnarray*}
so that we may summarize this as \begin{eqnarray}\label{cancel_u2-u1} xy^{u_1}xy^{u_2}xy^{a_i}y^{v_k}(xy^{u_t})\prod_{\substack {1\le j\le s\\j\ne k}}y^{v_j} = y^{\gamma+a_1-a_i-2v_k}\in \sigma(S)\end{eqnarray}

Let $A=\{\gamma+a_1-a_i+u_1-u_2:1\le i\le\ell\}$, and $B=\{-2v_1,\ldots,-2v_s\}$. Note that by lemma \ref{lm2_Sodd}, $|A|\ge \frac{t}{2}-1$; also,  $|B|\ge s/2$. Since $u_2-u_1\equiv\gamma\pmod 2$, we have $A, B\subseteq 2\cdot \Z_n$, and satisfy $|A|+|B|\ge\frac{n}{2}=|2\cdot\Z_n|$. 

 If $|A|+|B|>n/2$  then by lemma \ref{lm3} we must have $A+B=2\cdot \Z_n$, so that there exist $i, k$ such that $\gamma+a_1-a_i+u_1-u_2-2v_k=0$ in $\Z_n$. But this gives $\gamma+a_1-a_i-2v_k=u_2-u_1$. Hence we may  cancel  $xy^{u_1}xy^{u_2}=y^{u_2-u_1}$ from both sides of (\ref{cancel_u2-u1}) to get a one-product of length $n$ from the elements of $S$. That settles this case as well.

So, the only unsettled case corresponds to the one where $|A|+|B|=\frac{n}{2}$;  in particular, we must have $|A|=\frac{t}{2}-1$, and $|B|=\frac{s}{2}$. 

Note that by the same arguments as before, we also have 
$$xy^{u_1}xy^{u_2}xy^{a_i}y^{v_T}(xy^{u_t})\prod_{j\not\in T}y^{v_j}=y^{\gamma+a_1-a_i-2v_T}\in\sigma(S), \textrm{\ where\ }  v_T:=\sum_{i\in T} v_i\textrm{\ for\ }  T\subseteq[s].$$  In particular, if we set $\mathcal{S}(B):=\{2v_T: T\subseteq[s]\}$ then if $|A|+|\mathcal{S}(B)|>\frac{n}{2}$ also yields the same conclusion as above. Hence we may assume that $\mathcal{S}(B)=B$ but this implies that $B$ is a subgroup of $2\cdot\Z_n$.  
  
Finally,  since $|A|=\frac{t}{2}-1$, we have $$\left|\prod\nolimits_{t-3}(xy^{u_3},\ldots, xy^{u_{t-1}})\right| = \frac{t-2}{2}$$ so by lemma \ref{lm2_Sodd} we conclude that $\{u_3,\ldots,u_{t-1}\}$ is almost a coset of the subgroup of $\Z_n$ of order $t-2$. But since $t\ge\frac{n}{2}+2$, we must necessarily have $t=\frac{n}{2}+2$ and $s=\frac{n}{2}$. Consequently,  $B=2\cdot \Z_n$. Further, by lemma \ref{lm2_Sodd}, since $A$ is `almost'  a coset $K$ of $2\cdot\Z_n$, and since $A\subset 2\cdot \Z_n$, it follows that $A+B=2\cdot\Z_n$. But then clearly, $0\in A+B$, and the proof is complete as before.
 \end{itemize}\qed

\section{Concluding Remarks}
\begin{itemize}
\item The arguments in the proof of the main theorem may be slightly modified to also 
proves the following: For $n$ odd, any set  $S\subseteq D_{2n}$ of size $n+1$, admits a one-product subsequence of size $n$.
\item We believe that some of these ideas might extend to determine the Harborth constant of other semidirect products of cyclic groups.  An instructive first instance would be to consider the groups $C_p\rtimes \Z_n$ with $p$ prime and $p\mid\phi(n)$ where $\phi(n)$ denotes the Euler totient function. 
\item A more general invariant goes as follows. For a fixed integer $k$ and a finite group $G$, find the minimum integer $n$ such that every \textit{$k$-restricted sequence}, i.e., a sequence of elements from $G$ where no element appears more than $k$ times, of elements from $G$ admits a $1$-product subsequence of length $\exp(G)$. While this problem has been studied in the context of abelian groups, this appears to be an unstudied problem for general nonabelian groups. In particular, the aforementioned problem is open (to the best of our knowledge) even for dihedral groups, even for $k=2$. 

\end{itemize}
\section*{Acknowledgments}
We thank the reviewer for their careful reading and help with rewording and rephrasing many parts of the earlier draft, and also pointing to inconsistencies and inaccuracies that were present.

 \vspace{0.5cm}
 
 \newpage
 \section{Appendix: Remainder of the proof of lemma \ref{lm2_Sodd}}
 We furnish the remaining details of the proof of lemma \ref{lm2_Sodd}. We have already completed the proof in the case when the distinguished element of $A_0$ was a  generic former element. Now, we shall consider the remaining cases. 
\begin{itemize}
\item  The distinguished element is $\alpha_{s+2}-\alpha_{s+1}$ (the special former element corresponding to $i=1$):
Arguing as in the proof of the lemma, we get $\alpha_{s+3}-\alpha_2\in\{\alpha_{s+2}-\alpha_{s+1}+\frac{n}{2}, \alpha_{s+2}-\alpha_1\}$. If $\alpha_{s+3}-\alpha_2= \alpha_{s+2}-\alpha_1$, then a simple inductive argument shows that
$$\alpha_2-\alpha_1=\alpha_3-\alpha_2=\cdots=\alpha_{s+3}-\alpha_{s+2}= d, \textrm{\ say.}$$
By (\ref{alpha_3}) we have $\alpha_{s+3}-\alpha_{s+2}=\alpha_{s+2}-\alpha_s$ as well which gives $d=0$, a contradiction. Hence, $\alpha_{s+3}-\alpha_2=\alpha_{s+2}-\alpha_{s+1}+\frac{n}{2}$. Now, an inductive argument by forward propagation gives $\alpha_{s+\ell+1}-\alpha_2=\alpha_{s+\ell}-\alpha_1$ for all $3\le \ell\le s$. These along with (\ref{alpha_3}) give  $d=\alpha_{\ell +1}-\alpha_{\ell}$ for all $1\le\ell\le s-1$, and $\alpha_{s+3}-\alpha_{s+2}=\alpha_{s+2}-\alpha_s$. Thus, the upshot is 
$$\alpha_2-\alpha_1=\alpha_3-\alpha_2=\cdots=\alpha_s-\alpha_{s-1}=\alpha_{s+4}-\alpha_{s+3}=\cdots=\alpha_{2s+1}-\alpha_{2s}= d,$$ and $\alpha_{s+2}-\alpha_s=\alpha_{s+3}-\alpha_{s+2}$.

To complete the proof, we shall now show that $\alpha_{s+3}-\alpha_3 = \alpha_{s+2}-\alpha_1=\frac{n}{2}$. Once this is achieved, it is a straightforward check to see that 
$$\alpha_{s+1}-\alpha_s=\alpha_{s+2}-\alpha_{s+1}=d, \alpha_{s+3}-\alpha_{s+2}=2d$$ and that completes the analysis of this case.

To see why the claim holds, we observe that  $\alpha_{s+3}-\alpha_3\in\{\alpha_{s+2}-\alpha_1,\alpha_{s+2}-\alpha_2\}$. If $\alpha_{s+3}-\alpha_3=\alpha_{s+2}-\alpha_2 = d$, then in particular $\alpha_{s+2}-\alpha_1=sd=\frac{n}{2}$, and $\alpha_{s+2}-\alpha_s=d$ as well, so
$$ sd=\alpha_{s+3}-\alpha_2=\alpha_{s+2}-\alpha_{s+1}+\frac{n}{2}$$ which implies that $\alpha_{s+2}=\alpha_{s+1}$, a contradiction.

\item The distinguished element is $0$ (the special former element corresponding to $i=0$): In this case there are some slight differences in some of the details.  First, we have 
\begin{eqnarray*}
\alpha_{s+1}-\alpha_1 &=& \frac{n}{2}\\
\alpha_{s+2}-\alpha_{s+2-\ell}+\frac{n}{2} &=& \alpha_{s+1+\ell}-\alpha_1
\end{eqnarray*} 
which gives 
\begin{eqnarray}\label{alpha_4} \alpha_{s+2}-\alpha_{s+2-\ell}=\alpha_{s+\ell+1}-\alpha_{s+1}\textrm{\ for\ all\ }1\le\ell\le s.\end{eqnarray}
We claim that $\alpha_{s+3}-\alpha_2=\frac{n}{2}$. Indeed, it is again easy to see that $\alpha_{s+3}-\alpha_2\in\{\frac{n}{2},\alpha_{s+2}-\alpha_1\}$. If $\alpha_{s+3}-\alpha_2=\alpha_{s+2}-\alpha_1$, then a forward propagation inductive argument gives 
$$\alpha_2-\alpha_1=\alpha_{s+3}-\alpha_{s+2}=\cdots=\alpha_{2s+1}-\alpha_{2s}= d,\textrm{\ say}.$$
Combining this with (\ref{alpha_4}) gives 
\begin{eqnarray*} \alpha_{i+1}-\alpha_i=d\textrm{\ for\ all\ }1\le i\le s.\end{eqnarray*}
In particular, $sd=\frac{n}{2}$.

Finally, since $\alpha_{s+3}-\alpha_3\in\{\frac{n}{2}, \alpha_{s+2}-\alpha_2\}$, we observe that either possibility leads to a contradiction: If $\alpha_{s+3}-\alpha_3=\frac{n}{2}$ then $\alpha_{s+2}-\alpha_{s+1}+(s-1)d=\frac{n}{2}$ which forces  $\alpha_{s+2}-\alpha_{s+1}=d$ as well. But then $\alpha_{i+1}-\alpha_i=d$ for all $i$ and since $sd=\frac{n}{2}$, this gives $\alpha_{2s+1}=\alpha_1$, a contradiction. The other possibility leads to a similar contradiction. 

Hence, we conclude: $\alpha_{s+3}-\alpha_2=\frac{n}{2}$ as claimed.

This then leads via backward propagation to $\alpha_{s+3}-\alpha_{s+3-\ell}=\alpha_{s+2}-\alpha_{s+2-\ell}$ for all $1\le\ell\le s$ which in turn implies that 
$$ \alpha_3-\alpha_2=\alpha_4-\alpha_3=\cdots=\alpha_{s+3}-\alpha_{s+2}=d,\textrm{\ say}.$$ This observation in conjunction with (\ref{alpha_4}) gives 
$$ \alpha_{i+1}-\alpha_i = d\textrm{\ for\ all\ }2\le i\le 2s$$
which leaves only one gap undetermined, viz., $\alpha_2-\alpha_1$. But note that $(s+1)d=\alpha_{s+3}-\alpha_2=\frac{n}{2}$ which gives $d=\frac{n}{2(s+1)}$. Finally,
 \begin{eqnarray*}(s+1)d=\frac{n}{2}=\alpha_{s+1}-\alpha_1 &=& \alpha_2-\alpha_1 + \alpha_{s+1}-\alpha_2\\
  &=& \alpha_2-\alpha_1 +(s-1)d\end{eqnarray*} which gives $\alpha_2-\alpha_1=2d$.
  
\item The distinguished element is $\alpha_{s+2}-\alpha_2$ (the special latter element corresponding to $i=s$): It is easy to check that  
\begin{eqnarray*} \alpha_{s+2}-\alpha_{s+2-\ell}+\frac{n}{2}&=&\alpha_{s+2+\ell}-\alpha_1\textrm{\ for\ } 1\le\ell\le s-1.   \end{eqnarray*}         
Consequently, it is  trivial to see that both forward and backward propagation carry through all the way, so that we have 
$$\alpha_{i+1}-\alpha_i = d\textrm{\ for\  all\ } 1\le i\le 2s.$$ As $\alpha_{s+2}-\alpha_1 = \frac{n}{2}$, we then have that $(s+1)d=\frac{n}{2}$, and that settles this case.
 \item The distinguished element is $\alpha_{s+2}-\alpha_{3}$ (the special latter element corresponding to $i=s-1$): Here, backward propagation goes all the way yielding
 $$\alpha_2-\alpha_1=\cdots=\alpha_{s+3}-\alpha_{s+2}=d,\textrm{\ say}.$$ Forward propagation carries on till almost all the way giving
 $$\alpha_{s+2}-\alpha_{s+2-\ell}=\alpha_{s+2+\ell}-\alpha_{s+2}\textrm{ for\ }1\le\ell\le s-2$$
 which gives $$\alpha_{s+2}-\alpha_{s+1}=\cdots\alpha_{2s}-\alpha_{2s-1}=d.$$ In particular, $(s+1)d=\frac{n}{2}$.
 
 Finally, an argument similar to the ones in the preceding analyses (details omitted) gives $\alpha_{2s+1}-\alpha_2=\alpha_{s+2}-\alpha_3+\frac{n}{2}$ which implies that 
\begin{eqnarray*}\alpha_{2s+1}-\alpha_{2s} +(2s-2)d =  \alpha_{2s+1}-\alpha_2 &=& \alpha_{s+2}-\alpha_3+\frac{n}{2}\\
 &=& (s-1)d + (s+1)d.\end{eqnarray*} 
Consequently,  $\alpha_{2s+1}-\alpha_{2s}=2d$, and this case is hence settled. See figure 3 for an illustration of the specifics of each of these cases.

\end{itemize}
If the distinguished element is a latter element of $A_0$ (i.e., if the distinguished element is of the form $\alpha_{s+i}-\alpha_1$ for some $2\le i\le s+1$),  the argument is very similar. One can show that 
$$\left\{\alpha_1,\ldots,\alpha_{2s+1}\right\}=\{a,a+d,\ldots,a+(s+2-i)d,a+(s+4-i)d,\ldots, a+2sd\}$$ with $(s+1)d=\frac{n}{2}$, so that the missing element is $a+(s+3-i)d$. We omit the details.

\begin{figure}
\begin{pspicture}(0,-5.52)(15.817116,5.52)
\psline[linecolor=black, linewidth=0.04](0.19711548,5.28)(14.997115,5.28)(14.997115,5.28)
\psline[linecolor=black, linewidth=0.04](0.19711548,2.08)(14.997115,2.08)(14.997115,2.08)
\psline[linecolor=black, linewidth=0.04](0.19711548,-1.12)(14.997115,-1.12)
\psline[linecolor=black, linewidth=0.04](0.19711548,-4.32)(14.997115,-4.32)
\psdots[linecolor=black, dotsize=0.2](0.19711548,5.28)
\psdots[linecolor=black, dotsize=0.2](0.9971155,5.28)
\psdots[linecolor=black, dotsize=0.2](7.3971157,5.28)
\psdots[linecolor=black, dotsize=0.2](6.5971155,5.28)
\psdots[linecolor=black, dotsize=0.2](8.997115,5.28)
\psdots[linecolor=black, dotsize=0.2](9.797115,5.28)
\psdots[linecolor=black, dotsize=0.2](14.197116,5.28)
\psdots[linecolor=black, dotsize=0.2](14.997115,5.28)
\psdots[linecolor=black, dotsize=0.2](0.19711548,2.08)
\psdots[linecolor=black, dotsize=0.2](1.7971154,2.08)
\psdots[linecolor=black, dotsize=0.2](2.5971155,2.08)
\psdots[linecolor=black, dotsize=0.2](14.197116,2.08)
\psdots[linecolor=black, dotsize=0.2](14.997115,2.08)
\psdots[linecolor=black, dotsize=0.2](0.19711548,-1.12)
\psdots[linecolor=black, dotsize=0.2](0.9971155,-1.12)
\psdots[linecolor=black, dotsize=0.2](14.197116,-1.12)
\psdots[linecolor=black, dotsize=0.2](14.997115,-1.12)
\psdots[linecolor=black, dotsize=0.2](0.19711548,-4.32)
\psdots[linecolor=black, dotsize=0.2](0.9971155,-4.32)
\psdots[linecolor=black, dotsize=0.2](14.997115,-4.32)
\psdots[linecolor=black, dotsize=0.2](13.397116,-4.32)
\psdots[linecolor=black, dotsize=0.2](12.5971155,-4.32)
\rput[bl](0.59711546,5.4){$d$}
\rput[bl](6.9971156,5.4){$d$}
\rput[bl](8.197116,5.4){$2d$}
\rput[bl](9.397116,5.4){$d$}
\rput[bl](14.5971155,5.4){$d$}
\rput[bl](0.59711546,2.2){$2d$}
\rput[bl](2.1971154,2.2){$d$}
\rput[bl](14.5971155,2.2){$d$}
\rput[bl](0.59711546,-1.0){$d$}
\rput[bl](14.5971155,-1.0){$d$}
\rput[bl](0.59711546,-4.2){$d$}
\rput[bl](12.997115,-4.2){$d$}
\rput[bl](14.197116,-4.2){$2d$}
\rput[bl](0.19711548,4.8){$\alpha_1$}
\rput[bl](0.9971155,4.8){$\alpha_2$}
\rput[bl](6.5971155,4.8){$\alpha_{s+1}$}
\rput[bl](7.3971157,4.8){$\alpha_{s+2}$}
\rput[bl](8.997115,4.8){$\alpha_{s+3}$}
\rput[bl](9.797115,4.8){$\alpha_{s+4}$}
\rput[bl](14.197116,4.8){$\alpha_{2s}$}
\rput[bl](14.997115,4.8){$\alpha_{2s+1}$}
\rput[bl](0.19711548,1.6){$\alpha_1$}
\rput[bl](1.7971154,1.6){$\alpha_2$}
\rput[bl](2.5971155,1.6){$\alpha_3$}
\rput[bl](14.197116,1.6){$\alpha_{2s}$}
\rput[bl](14.997115,1.6){$\alpha_{2s+1}$}
\rput[bl](0.19711548,-1.6){$\alpha_1$}
\rput[bl](0.9971155,-1.6){$\alpha_2$}
\rput[bl](14.197116,-1.6){$\alpha_{2s}$}
\rput[bl](14.997115,-1.6){$\alpha_{2s+1}$}
\rput[bl](0.19711548,-4.8){$\alpha_1$}
\rput[bl](0.9971155,-4.8){$\alpha_2$}
\rput[bl](12.197116,-4.8){$\alpha_{2s-1}$}
\rput[bl](13.397116,-4.8){$\alpha_{2s}$}
\rput[bl](14.997115,-4.8){$\alpha_{2s+1}$}
\rput[bl](3.3971155,3.68){Distinguished element: $\alpha_{s+2}-\alpha_{s+1}$.}
\rput[bl](3.3971155,0.88){Distinguished element: $0$.}
\rput[bl](3.3971155,-2.32){Distinguished element: $\alpha_{s+2}-\alpha_2$.}
\rput[bl](3.3971155,-5.52){Distinguished element: $\alpha_{s+3}-\alpha_2$.}
\rput[bl](9.397116,3.68){(Special former element for $i=1$.)}
\rput[bl](9.397116,0.88){(Special former element for $i=0$.)}
\rput[bl](9.397116,-2.32){(Special latter element for $i=s$.)}
\rput[bl](9.397116,-5.52){(Special latter element for $i=s-1$.)}
\end{pspicture}
\vspace{1cm}

\caption{The $\alpha_i$ in order for the case where the distinguished element is a special element of $A_0$.}
\end{figure}

\end{document}